\documentclass[12pt]{article}

\usepackage{amsfonts,amsmath,amssymb,amsthm}
\usepackage{epsfig}
\usepackage{graphics}
\usepackage{graphicx}
\usepackage{times}

\title{The zeros of the derivative of the Riemann zeta function near the
critical line}
\begin{document}
\newcommand\tabcaption{\def\@captype{table}\caption}
\newtheorem{sthrm}{Theorem}[section]
\newtheorem{thrm}{Theorem}
\newtheorem{athrm}{Theorem}
\newtheorem{cor}[thrm]{Corollary}
\newtheorem{co}{Corollary}
\newtheorem{ex}{Example}
\newtheorem{lemma}{Lemma}[section]
\newtheorem{prop}{Proposition}[section]
\newtheorem*{thrma}{Theorem A}
\newtheorem*{thrmb}{Theorem B}
\newtheorem*{thrmc}{Theorem C}
\newtheorem*{thrmd}{Theorem D}
\newtheorem*{thrme}{Theorem E}
\newtheorem*{claim}{Claim}
\newtheorem*{con}{Conjecture}
\newtheorem*{thrmpl}{Phragm\'en-Lindel\"of Theorem}

\numberwithin{equation}{section}

\def\R{\text{\rm{Re}}\,}
\def\I{\text{\rm{Im}}\,}
\def\s{\text{\rm sign}}
\def\bysame{\leavevmode\hbox to3em{\hrulefill}\thinspace}

\author{Haseo Ki\\
\\
Department of Mathematics, Yonsei
University, Seoul 120--749, Korea\\
\\
{\it haseo\@yonsei.ac.kr}}

\parskip=12pt

\maketitle 

\begin{abstract}
We study the horizontal distribution of zeros of $\zeta'(s)$ which
are denoted as $\rho'=\beta'+i\gamma'$. We assume the Riemann
hypothesis which implies $\beta'\geqslant1/2$ for any non-real zero
$\rho'$, equality being possible only at a multiple zero of
$\zeta(s)$. In this paper we prove that
$\liminf\left(\beta'-1/2\right)\log\gamma'\not=0$ if and only if for
any $c>0$ and $s=\sigma+it$ with $|\sigma-1/2|<c/\log t$
$(t\geqslant10)$
$$
\frac{\zeta'}{\zeta}(s)=\frac{1}{s-\rho}+O(\log t),
$$
where $\rho=1/2+i\gamma$ is the closest zero of $\zeta(s)$ to $s$
and the origin. We also show that if
$\liminf\left(\beta'-1/2\right)\log\gamma'\not=0$, then for any
$c>0$ and $s=\sigma+it$ ($t\geqslant10$), we have
$$
\log\zeta(s)=O\left(\frac{(\log t)^{2-2\sigma}}{\log\log t}\right)
$$
uniformly for $1/2+c/\log t\leqslant\sigma\leqslant\sigma_1<1$.
\end{abstract}


\pagestyle{myheadings}

\section{Introduction}
The Riemann hypothesis (RH) states that the real part of any nonreal
zero of the Riemann zeta function $\zeta(s)$ is $1/2$. A. Speiser
\cite{Sp} has a theorem which says that RH is equivalent to the
nonexistence of nonreal zeros of $\zeta'(s)$ in $\R(s)<1/2$.

We let $\rho'=\beta'+i\gamma'$ denote a zero of $\zeta'(s)$ where a
sum over $\rho'$ is repeated according to multiplicity. N. Levinson
and H. L. Montgomery \cite{LM} proved that for nonreal zeros of
$\zeta'(s)$ the average value of $\beta'$ with $0<\gamma'\leqslant
T$ is $1/2+\log\log T/\log T$. But it is likely that $\beta'-1/2$ is
usually of the order $1/\log T$ rather than $\log\log T/\log T$.

By a result of B. C. Berndt \cite{Be}, the number of zeros with
ordinate less than $T$ is
$$
\sum_{0<\gamma'\leqslant T}1=\frac{T}{2\pi}\log\frac{T}{4\pi
e}+O(\log T).
$$
Under RH, K. Soundararajan \cite{So} demonstrated the presence of a
positive proportion of zeros of $\zeta'(s)$ in the region
$\sigma<1/2+\nu/\log T$ for all $\nu\geqslant2.6$. `A positive
proportion of zeros' means
$$
\liminf_{T\to\infty}\frac{1}{\frac{T}{2\pi}\log T}
\#\{\rho':\beta'\leqslant1/2+\nu/\log T,0<\gamma'\leqslant T\}>0.
$$
In his remarkable work, Y. Zhang \cite{Zh} showed that not only
unconditionally there exists a $\nu>0$ such that a positive
proportion of zeros of $\zeta'(s)$ are in the region
$|\sigma-1/2|<\nu/\log T$, but also assuming RH and a strong
hypothesis of the distribution of zeros of $\zeta(s)$, $\nu$ can be
arbitrary small. Recently, Feng \cite{Fe} proved the Zhang's second
theorem only assuming the strong hypothesis. Thus it is very
probable that
$$
\liminf\left(\beta'-1/2\right)\log\gamma'=0.
$$

Assuming the truth of RH, K. Soundararajan \cite{So} conjectured
that the following two statements are equivalent:

{\rm (i)} $\liminf(\beta'-1/2)\log\gamma'=0$;

{\rm (ii)} $\liminf(\gamma^+-\gamma)\log\gamma=0$ where $\gamma^+$
is the least ordinate of a zero of $\zeta(s)$ with
$\gamma^+>\gamma$.

Y. Zhang \cite{Zh} has shown that (ii) implies (i) as follows.

\begin{thrma}Assume RH. Let $\alpha_1$ and $\alpha_2$ be
positive constants satisfying $\alpha_1<2\pi$ and
$$\alpha_2>\alpha_1\left(1-\sqrt{\frac{\alpha_1}{2\pi}}\right)^{-1}.$$
If $\rho=1/2+i\gamma$ is a zero of $\zeta(s)$ such that $\gamma$ is
sufficiently large and $\gamma^+-\gamma<\alpha_1(\log\gamma)^{-1}$,
then there exists a zero $\rho'$ of $\zeta'(s)$ such that
$$|\rho'-\rho|<\alpha_2(\log\gamma)^{-1}.$$
\end{thrma}

In this paper, we consider the converse of Theorem A. Namely is it
true that (i) implies (ii)? Concerning this problem, we have the
following.

\begin{thrm}\label{t:main} Assume RH. Then the following are equivalent:

{\rm (1)} $\liminf(\beta'-1/2)\log\gamma'\not=0$;

{\rm (2)} For any $c>0$ and $s=\sigma+it$ with $|\sigma-1/2|<c/\log
t$ $(t\geqslant10)$
$$\frac{\zeta'}{\zeta}(s)=\frac{1}{s-\rho}+O(\log t),$$
where $\rho=1/2+i\gamma$ is the closest zero of $\zeta(s)$ to $s$
and the origin;
\end{thrm}

\begin{co}\label{co:1} Assume RH and
$\liminf(\beta'-1/2)\log\gamma'\not=0$. Then, for any $c>0$ and
$s=\sigma+it$ $(t\geqslant10)$, we have
$$
\frac{\zeta'}{\zeta}(s)=O\left((\log t)^{2-2\sigma}\right)
$$
uniformly for $1/2+c/\log t\leqslant\sigma\leqslant\sigma_1<1$.
\end{co}

Based on our theorems and Soundararajan's conjecture, we speculate
as follows.

\begin{con}\label{con:1} Assume RH. Then the following are equivalent:

{\rm (i)${}'$} For any $c>0$ and $s=\sigma+it$ $(t\geqslant10)$, we
have
$$
\frac{\zeta'}{\zeta}(s)=O\left((\log t)^{2-2\sigma}\right)
$$
uniformly for $1/2+c/\log t\leqslant\sigma\leqslant\sigma_1<1$;

{\rm (ii)${}'$} The negation of {\rm (ii)}, i.e.,
$\liminf(\gamma^+-\gamma)\log\gamma\not=0$.
\end{con}

We briefly introduce why (2) in Theorem \ref{t:main} doesn't seem
possible. We let $s=\sigma+it$ for real numbers, $\sigma, t$. It is
known in \cite[p. 99]{Da} and \cite[Theorem 9.6(A)]{Ti} that for
$t\geqslant10$ and $-1\leqslant\sigma\leqslant2$ we have
\begin{equation}
\frac{\zeta'}{\zeta}(s)=\sum_{|\gamma-t|\leqslant1}\frac{1}{s-\rho}+O(\log
t),
\end{equation}
where the sum is over the ordinates of the complex zeros
$\rho=\beta+i\gamma$ of $\zeta(s)$. However, assuming RH, we may
expect a stronger result. Under RH we have
\begin{equation}
\frac{\zeta'}{\zeta}(s)=\sum_{|\gamma-t|\leqslant1/\log\log
t}\frac{1}{s-\rho}+O(\log t)
\end{equation} for $t\geqslant10$.
For this we refer to \cite[p. 357 (14.15.2)]{Ti}. According to the
formula (1.2), we can see that the formula (2) in Theorem 1 is very
unrealistic, because `$|\sigma-1/2|<c/\log t$' in (2) doesn't seem
probable instead of `$|\gamma-t|\leqslant1/\log\log t$' in (1.2).

Concerning Corollary 1, it is worth noting that under RH
\begin{equation}
\frac{\zeta'}{\zeta}(s)=O\left((\log t)^{2-2\sigma}\right)
\end{equation}
holds uniformly for $1/2+c/\log\log
t\leqslant\sigma\leqslant\sigma_1<1$, where $t\geqslant10$, $c>0$
and `$O$' depends upon $c$ and $\sigma_1$. For the proof of it, we
apply the fact \cite[(14.14.5)]{Ti} to the formula
$$
\frac{\zeta'(s)}{\zeta(s)}=\frac{1}{2\pi
i}\int_{|s-z|=\frac{c}{2\log\log t}}\frac{\log\zeta(z)}{(s-z)^2}dz.
$$ Then we
immediately get (1.3). In proving (1.3), one cannot relax the
condition `$1/2+c/\log\log t$' as
$$
\frac{1}{2}+\frac{c}{\log t}\qquad(t\to\infty).
$$
With this information, it is very likely that on $\R(s)=1/2+c/\log
t$  $(t\geqslant10)$,
$$
\frac{\zeta'}{\zeta}(s)\not=O(\log t).
$$

On the other hand, Corollary 1 follows from a Phragm\'en-Lindel\"of
argument, provided that we have
$$
\frac{\zeta'}{\zeta}(s)=O(\log t)
$$
on $\R(s)=1/2+c/\log t$  $(t\geqslant10)$. Thus we find the behavior
of the logarithmic derivative of the Riemann zeta function near
$\R(s)=1/2$ subtle and so we need a deep observation about the
Riemann zeta function near the critical line to establish
$\liminf(\beta'-1/2)\log\gamma'=0$. In fact we will see from Theorem
4 in Section 2 that the behavior of $\zeta'(s)/\zeta(s)$ on
$\R(s)=1/2+c/\log t$ is very much related to
$$
\sum_{0<|\gamma-\tilde{\gamma}|<1}\frac{1}{\gamma-\tilde{\gamma}},
$$
where $1/2+i\gamma,1/2+i\tilde{\gamma}$ are complex zeros of
$\zeta(s)$ and the sum is over $\tilde{\gamma}$.

From Corollary 1, we can demonstrate the following.

\begin{co}\label{co:2} Assume RH and
$\liminf(\beta'-1/2)\log\gamma'\not=0$. Then, For any $c>0$ and
$s=\sigma+it$ $(t\geqslant10)$, we have
$$
\log\zeta(s)=O\left(\frac{(\log t)^{2-2\sigma}}{\log\log t}\right)
$$
uniformly for $1/2+c/\log t\leqslant\sigma\leqslant\sigma_1<1$.
\end{co}

Assuming RH, it is known in \cite[p. 355 (14.14.5)]{Ti} and
\cite[Theorem 14.14(B)]{Ti} that for any $c>0$ and $s=\sigma+it$
$(t\geqslant10)$, we have
\begin{equation}
\log\zeta(s)=O\left(\frac{(\log t)^{2-2\sigma}}{\log\log t}\right)
\end{equation}
holds uniformly for $1/2+c/\log\log
t\leqslant\sigma\leqslant\sigma_1<1$, and there exists an absolute
constant $c^*>0$ depending on $c$ such that
\begin{equation}
-c^*\frac{\log t}{\log\log
t}\log\left(\frac{2}{\left(\sigma-\frac{1}{2}\right)\log\log
t}\right)<\log|\zeta(s)|<c^*\frac{\log t}{\log\log t}
\end{equation}
holds for $1/2<\sigma\leqslant1/2+c/\log\log t$;
\begin{equation}
\arg\zeta(s)=O\left(\frac{\log t}{\log\log t}\right)
\end{equation}
holds uniformly for $1/2\leqslant\sigma\leqslant1/2+c/\log\log t$.
However, as in the proof of (1.3), one cannot relax the condition
`$1/2+c/\log\log t$' of (1.4) as $1/2+c/\log t\,\ (t\to\infty)$. On
the other hand, we have $\Omega$-theorems related to $\log\zeta(s)$
near the critical line. For these, we refer to Montgomery's results
\cite{Mo2} and \cite[p. 209]{Ti}. In particular, Montgomery showed
that assuming RH, for $1/2\leqslant\sigma<1$ and any real $\theta$,
there is a $t$ with $T^{1/6}\leqslant t\leqslant T$ such that
$$
\R\left(e^{-i\theta}\log\zeta(s)\right)\geqslant\frac{1}{20}(\log
T)^{1-\sigma}(\log\log T)^{-\sigma}.
$$
See \cite[p. 512]{Mo2}. Recently, from random matrix theory, in the
case that $\theta=0$ in the above $\Omega$-result, it is conjectured
in \cite{FGH} that we have the following:
$$
\max_{t\in[0,T]}\left|\zeta\left(\frac{1}{2}+it\right)\right|=\exp\left((1+o(1))\sqrt{\frac{1}{2}\log
T\log\log T}\right).
$$
Concerning negative values of $\log|\zeta(s)|$, we observe that by
(1.5) and (1.6), it is possible that we have
$$
\left|\log\zeta\left(\frac{1}{2}+\frac{1}{\log
t}+it\right)\right|>\psi(t)\frac{\log t}{\log\log t}
$$
for some arbitrarily large values of $t$, where $\psi(t)\to\infty$
as $t\to\infty$. A sharp $\Omega$-result like this implies
$\liminf(\beta'-1/2)\log\gamma'=0$. Namely, we note that it will
follow if we show that (1.4) does not hold uniformly for $1/2+1/\log
t\leqslant\sigma\leqslant1/2+1/\log\log t$. Thus, Corollary 2 is
useful in investigating the horizontal behavior of zeros of the
derivative of the Riemann zeta function.

We apply our theorem to mean values of the logarithmic derivative of
the Riemann zeta function:
$$
\int_0^T\left|\frac{\zeta'}{\zeta}(\sigma+it)\right|^2dt
$$
for $\sigma=1/2+a/\log T$.

We may get the following from a result of A. Selberg \cite[equation
(1.2)]{Se}.

\begin{thrmb}\label{t:selberg}
Assume RH. Then
$$\int_0^T\left|\frac{\zeta'}{\zeta}(\sigma+it)\right|^2dt\sim\frac{1}{4a^2}T\log^2T$$
holds where $\sigma=1/2+a/\log T$ and $a\to\infty$, $a=o(\log T)$.
\end{thrmb}

We introduce more studies on mean values of the logarithmic
derivative of $\zeta(s)$. Following Montgomery \cite{Mo1}, let
$$
F(\alpha,T)=\frac{1}{\frac{T}{2\pi}\log
T}\sum_{0<\gamma,\tilde{\gamma}\leqslant
T}T^{i\alpha(\gamma-\tilde{\gamma})} w(\gamma-\tilde{\gamma}),
$$
where $\beta+i\gamma$ and $\tilde{\beta}+i\tilde{\gamma}$ are zeros
of $\zeta(s)$ and $w(u)=\frac{4}{4+u^2}$. Montgomery conjectured
that for any fixed $A>1$,
\begin{displaymath}
F(\alpha,T)\sim1\qquad\text{ uniformly for
}1\leqslant\alpha\leqslant A.\leqno({\rm MH})
\end{displaymath}

Under RH, D. A. Goldston, S. M. Gonek and H. L. Montgomery \cite{Go}
considered
$$\int_0^T\left|\frac{\zeta'}{\zeta}(\sigma+it)\right|^2dt$$
for $\sigma=1/2+a/\log T$ as $a\to0$ and proved the following,
provided that (MH) is valid.

\begin{thrmc}\label{t:ggm} Assume RH and (MH). Then
$$
\int_0^T\left|\frac{\zeta'}{\zeta}(\sigma+it)\right|^2dt\sim\frac{1}{2a}T\log^2T
$$
holds where $\sigma=1/2+a/\log T$ and $a=a(T)\to0$ (sufficiently
slowly) as $T\to\infty$.
\end{thrmc}

We have the following as in Theorem \ref{t:ggm}.

\begin{thrm} Assume RH and
$\liminf(\beta'-1/2)\log\gamma'\not=0$. Then
$$
\int_0^T\left|\frac{\zeta'}{\zeta}(\sigma+it)\right|^2dt
=\frac{1}{2a}T\log^2T\left(1-\frac{\log(2\pi e)}{\log T}
+O\left(a\log\frac{1}{a}\right)+O\left(\frac{1}{T}\right) \right)
$$
holds where $\sigma=1/2+a/\log T$ and $a\to0$. Here `$O$' doesn't
depend upon $a$ and $T$.
\end{thrm}

(MH) implies Montgomery's pair correlation conjecture. That is,
$$\aligned
\sum_{\substack{0<\gamma,\tilde{\gamma}\leqslant T\\
0<\gamma-\tilde{\gamma}\leqslant2\pi\beta/\log
T}}1\sim\frac{T}{2\pi}\log T\int_0^{\beta} 1-\left(\frac{\sin\pi
u}{\pi u}\right)^2du.
\endaligned$$ Clearly (MH) implies
$$
\liminf(\gamma^+-\gamma)\log\gamma=0.
$$
Assuming RH, Theorem A says that
$$
\liminf(\beta'-1/2)\log\gamma'\not=0\text{ implies }
\liminf(\gamma^+-\gamma)\log\gamma\not=0.
$$
Thus the assumptions of Theorem C and Theorem 2 are contradictory to
each other. However, we have the similar conclusion in Theorem C and
Theorem 2. Further, a theorem of D. A. Goldston, S. M. Gonek and H.
L. Montgomery \cite[Theorem 3]{Go} says the following.

\begin{thrmd} Assume RH. Then, for $\sigma=1/2+a/\log T$ and
any fixed $a>0$,
$$
\int_0^T\left|\frac{\zeta'}{\zeta}(\sigma+it)\right|^2dt\sim\frac{1-e^{-2a}}{4a^2}T\log^2T
$$
holds as $T\to\infty$ if and only if the pair correlation conjecture
is true.
\end{thrmd}

Combining Theorem B and Theorem 2, we immediately have the
following.

\begin{thrm} \label{t:test}Assume RH and
$\liminf(\beta'-1/2)\log\gamma'\not=0$. Then, for $\sigma=1/2+a/\log
T$,
$$
\int_0^T\left|\frac{\zeta'}{\zeta}(\sigma+it)\right|^2dt\sim\frac{1-e^{-2a}}{4a^2}T\log^2T
$$
holds where $a\to0$ or $a\to\infty$ and $a=o(\log T)$.
\end{thrm}

Apparently, Theorem D and Theorem 3 are similar. However the
conclusion of Theorem D says a much stronger statement than that of
Theorem 3. We notice that the pair correlation conjecture implies
$\liminf(\beta'-1/2)\log\gamma'=0$. Thus the conclusion of Theorem 3
under $\liminf(\beta'-1/2)\log\gamma'\not=0$ is instructive to
understand the behavior of the mean value of the logarithmic
derivative of the Riemann zeta function.

\section{Proof of Theorem 1 and Corollaries 1, 2}

We arrange the zeros of the Riemann zeta function $\zeta(s)$ on the
upper half-plane as $\rho_1,\rho_2,\ldots$ with
$\rho_n=\beta_n+i\gamma_n$ and
$$
0<\gamma_1\leqslant\gamma_2\leqslant\ldots,
$$
where it appears precisely $m$ times consecutively in the above
sequence, if a zero is multiple with multiplicity $m$. RH is that
$\beta_n=1/2$ for any $n=1,2,3,\ldots$.

We state basic facts.

\begin{prop} Let $T>0$. Then, we have:

{\rm (1)} The number of zeros of $\zeta(s)$ in $0<\I(s)\leqslant T$
is
$$
\sum_{1<\gamma_n\leqslant T}1=\frac{T}{2\pi}\log\frac{T}{2\pi
e}+O(\log T);
$$

{\rm (2)} The number of zeros of $\zeta(s)$ in
$T\leqslant\I(s)\leqslant T+1$ is $O(\log T)$.
\end{prop}

For Proposition 2.1(1), see \cite[p. 98]{Da} and \cite[Theorem
9.4]{Ti}. Proposition 2.1(2) immediately follows from (1).

We start with the following theorem.

\begin{thrm}\label{t:lemma} Assume RH and
$\liminf(\gamma_{n+1}-\gamma_n)\log\gamma_n>0$. Then the following
three statements are equivalent:

{\rm (A)} $\liminf(\beta'-1/2)\log\gamma'>0$;

{\rm (B)} Let $c>0$ and $s=\sigma+it$. For sufficiently large $n$,
we have
$$
\frac{\zeta'}{\zeta}(s)=\frac{1}{s-\rho_n}+O(\log t),
$$ where
$\frac{\gamma_{n-1}+\gamma_n}{2}<t\leqslant\frac{\gamma_{n+1}+\gamma_n}{2}$
and $|\sigma-1/2|<c/\log\gamma_n$;

{\rm (C)} $\limsup |M_n|/\log\gamma_n<\infty$, where
$$
M_n=\sum_{0<|\gamma_m-\gamma_n|\leqslant1}
\frac{1}{\gamma_n-\gamma_m}.
$$
\end{thrm}

\begin{proof}[Proof of Theorem \ref{t:lemma}] We may assume that all but finitely many
nontrivial zeros of $\zeta(s)$ are simple and on $\R(s)=1/2$,
because we assume RH and
$\liminf(\gamma_{n+1}-\gamma_n)\log\gamma_n>0$.
\medskip

(C) $\Rightarrow$ (B). We recall that
$$
\frac{\zeta'}{\zeta}(s)=O(\log
t)+\sum_{|\gamma-t|\leqslant1}\frac{1}{s-\rho}
$$
for $-1\leqslant\R(s)\leqslant2$ and $t\geqslant2$.

\begin{prop} \label{p:2}Let $\delta>0$. Suppose that
$\liminf(\gamma_{n+1}-\gamma_n)\log\gamma_n>\delta$. Then we have
$$
\sum_{m\not=n}
\frac{1}{(\gamma_n-\gamma_m)^2}=O\left(\frac{\log^2\gamma_n}{\delta^2}\right)
$$
for sufficiently large $n$.
\end{prop}

\begin{proof}[Proof of Proposition 2.2] We write
\begin{equation}\begin{split}
\sum_{m\not=n}\frac{1}{(\gamma_m-\gamma_n)^2}=&\sum_{|\gamma_m-\gamma_n|>1}
\frac{1}{(\gamma_m-\gamma_n)^2}+
\sum_{0<|\gamma_m-\gamma_n|\leqslant1}\frac{1}{(\gamma_m-\gamma_n)^2}\\
=&I+II.
\end{split}\end{equation}
By Proposition 2.1(2), the number of $\gamma_n$'s between $t$ and
$t+1$ is $O(\log t)$ for $t>1$. Then, for some $C>0$, we get
\begin{equation}\begin{split}
I=&\sum_{k=1}^{\infty}\sum_{\gamma_n+k<\gamma_m\leqslant\gamma_n+k+1}\frac{1}{(\gamma_m-\gamma_n)^2}
+\sum_{k=1}^{\infty}\sum_{\gamma_n-k-1\leqslant\gamma_m<\gamma_n-k}\frac{1}{(\gamma_m-\gamma_n)^2}\\
\leqslant&\sum_{k=1}^{\infty}\frac{C\log(\gamma_n+k)}{k^2}+\sum_{\gamma_n-k>1}\frac{C\log(\gamma_n-k)}{k^2}\\
\leqslant&C\sum_{k=1}^{\infty}\frac{\log\gamma_n+\log
k}{k^2}+C\sum_{k=1}^{\infty}\frac{\log\gamma_n}{k^2}=O(\log\gamma_n)
\end{split}\end{equation}
By the assumption of Proposition 2.2, there exists a positive
integer $n_1$ such that
$$
\gamma_{m+1}-\gamma_m\geqslant\frac{\delta}{\log\gamma_m}\qquad\text{for
all }m\geqslant n_1.
$$
Using this, for sufficiently large $n$, we have
$$|\gamma_{n+k}-\gamma_n|\geqslant\frac{|k|\delta}{\log(\gamma_n+1)}\geqslant\frac{|k|\delta}{2\log\gamma_n}$$
for $|\gamma_{n+k}-\gamma_n|\leqslant1$ and $|k|\geqslant1$. Here
there exists a $a>0$ such that $1\leqslant|k|\leqslant
a\log\gamma_n$. Thus we get
$$
II\leqslant\sum_{1\leqslant|k|\leqslant a\log\gamma_n}
\frac{1}{\left(\frac{k\delta}{2\log\gamma_n}\right)^2}<\sum_{k=1}^{\infty}
\frac{2}{\left(\frac{k\delta}{2\log\gamma_n}\right)^2}=O\left(\frac{\log^2\gamma_n}{\delta^2}\right).
$$
We apply this inequality and (2.2) to (2.1) and then we obtain
$$\sum_{m\not=n}\frac{1}{(\gamma_m-\gamma_n)^2}=O\left(\frac{\log^2\gamma_n}{\delta^2}\right)$$
for sufficiently large $n$. This proves Proposition 2.2.
\end{proof}

\begin{prop} Choose $c_1>0$ such that
$\liminf(\gamma_{n+1}-\gamma_n)\log\gamma_n>2c_1$. Define $M_n(t)$
by
$$M_n(t)=\sum_{0<|\gamma_m-\gamma_n|\leqslant1}\frac{1}{t-\gamma_m}$$
for $\gamma_{n-1}+c_1/\log\gamma_n\leqslant
t\leqslant\gamma_{n+1}-c_1/\log\gamma_n$. Then we have
$$M_n(t)=O(\log t),$$
where the implied constant is absolute.
\end{prop}

\begin{proof}[Proof of Proposition 2.3] Since
$$M_n'(t)=\sum_{0<|\gamma_m-\gamma_n|<1}\frac{-1}{(t-\gamma_m)^2}<0,$$
$M_n(t)$ is decreasing. Thus it suffices to consider the endpoints
for the proof. By Proposition 2.2 and our assumption (C), we have
$$\aligned
M_n\left(\gamma_{n+1}-\frac{c_1}{\log\gamma_n}\right)=
&\sum_{0<|\gamma_m-\gamma_n|\leqslant1}\frac{1}{\gamma_{n+1}-\gamma_m-\frac{c_1}{\log\gamma_n}}\\
=&\sum_{0<|\gamma_m-\gamma_{n+1}|\leqslant1}\frac{1}{\gamma_{n+1}-\gamma_m-\frac{c_1}{\log\gamma_n}}+O(\log\gamma_n)\\
=&\sum_{0<|\gamma_m-\gamma_{n+1}|\leqslant1}\frac{1}{\gamma_{n+1}-\gamma_m-\frac{c_1}{\log\gamma_n}}
-\frac{1}{\gamma_{n+1}-\gamma_m}+\\
&M_{n+1}+O(\log\gamma_n)\\
=&O\left(\frac{1}{\log\gamma_n}\sum_{0<|\gamma_m-\gamma_{n+1}|<1}\frac{1}{(\gamma_{n+1}-\gamma_m)^2}\right)+O(\log\gamma_{n+1})\\
=&O(\log\gamma_n).\endaligned$$ Similarly, we have
$$M_n\left(\gamma_{n-1}+\frac{c_1}{\log\gamma_n}\right)=O(\log\gamma_n).$$
Proposition 2.3 follows.
\end{proof}

Using Proposition 2.1(2), Proposition 2.3 and (C), for
$\rho_n=1/2+i\gamma_n$ and $s=\sigma+it$ $(t\geqslant10)$, we get
$$\aligned
\frac{\zeta'}{\zeta}(s)=&\frac{1}{s-\rho_n}+\sum_{\substack{\rho\neq\rho_n\\0<|\gamma-t|\leqslant1}}\frac{1}{s-\rho}+O(\log t)\\
=&\frac{1}{s-\rho_n}+\sum_{0<|\gamma_m-\gamma_n|\leqslant1}\frac{1}{s-\rho_m}+
\sum_{\substack{|\gamma_m-\gamma_n|>1\\0<|\gamma_m-t|\leqslant1}}\frac{1}{s-\rho_m}+O(\log t)\\
=&\frac{1}{s-\rho_n}+\sum_{0<|\gamma_m-\gamma_n|\leqslant1}\frac{1}{s-\rho_m}+
O\left(\sum_{0<|\gamma-t|\leqslant1}1\right)+O(\log t)\\
=&\frac{1}{s-\rho_n}+\sum_{0<|\gamma_m-\gamma_n|\leqslant1}\frac{1}{s-\rho_m}-\frac{1}{i(t-\gamma_m)}+\frac{M_n(t)}{i}+O(\log t)\\
=&\frac{1}{s-\rho_n}+O\left(\sum_{m\not=n}\frac{\sigma-\frac{1}{2}}{(\gamma_m-\gamma_n)^2}\right)+O(\log
t)\endaligned$$ in $|\sigma-1/2|<c/\log\gamma_n.$ By this and
Proposition 2.2, (B) follows.
\medskip

(B) $\Rightarrow$ (A). We need the following proposition for this.

\begin{prop} Assume RH and
$\liminf(\gamma_{n+1}-\gamma_n)\log\gamma_n>0$. Let $\tilde{\delta}$
be such that $0<\tilde{\delta}<\delta$ where $\delta$ is as in
Proposition 2.2. Suppose $\zeta'(\beta'+i\gamma')=0$ for
sufficiently large $\gamma'$. If $|\gamma'-\gamma_n|\geqslant
\tilde{\delta}/(2\log\gamma')$ for all $n$, then we have
$$\frac{1}{2}\log\gamma'+O(1)=O\left(\left(\beta'-\frac{1}{2}\right)\frac{\log^2\gamma'}{\tilde{\delta}^2}\right).$$
\end{prop}

\begin{proof}[Proof of Proposition 2.4] We set
$$\xi(s)=\frac{s(s-1)}{2}\pi^{-s/2}\Gamma(s/2)\zeta(s)$$
for any $s\in\mathbb{C}$. It is known that for some constants $A$
and $B$,
$$\xi(s)=e^{A+Bs}\prod_{\rho}\left(1-\frac{s}{\rho}\right)e^{\frac{s}{\rho}},$$
where $\rho$ runs through all zeros of $\xi(s)$. See \cite[p.
80]{Da} and \cite[p. 30]{Ti} for this. We note that any zero $\rho$
of $\xi(s)$ is either $1/2+ i\gamma_n$ or $1/2-i\gamma_n$ for some
$n$, provided that RH is true. By the product formula of $\xi(s)$,
we get
$$\frac{\zeta'(s)}{\zeta(s)}=\frac{1}{2}\log\pi-\frac{1}{s-1}-\frac{1}{2}\frac{\Gamma'\left(\frac{s}{2}+1\right)}
{\Gamma\left(\frac{s}{2}+1\right)}+B+\sum_{\rho}\frac{1}{s-\rho}+\frac{1}{\rho}.$$
Then we have
\begin{equation}
\R\frac{\zeta'(s)}{\zeta(s)}=\frac{1}{2}\log\pi-\R\frac{1}{s-1}-\frac{1}{2}\R\frac{\Gamma'\left(\frac{s}{2}+1\right)}
{\Gamma\left(\frac{s}{2}+1\right)}+\sum_{\rho}\R\frac{1}{s-\rho}.
\end{equation}
We assumed that for any $n$,
\begin{equation}
|\gamma_n-\gamma'|>\frac{\tilde{\delta}}{2\log\gamma'}.
\end{equation}
By (2.3) we obtain that for $s=\beta'+i\gamma'$, we get
\begin{equation}
\sum_{\rho}\R\frac{1}{s-\rho}=-\frac{1}{2}\log\pi+\R\frac{1}{s-1}+\frac{1}{2}\R\frac{\Gamma'\left(\frac{s}{2}+1\right)}
{\Gamma\left(\frac{s}{2}+1\right)}.
\end{equation}
Applying the standard fact \cite[p. 73]{Da}
$$\frac{\Gamma'(s)} {\Gamma(s)}=\log
s+O(\frac{1}{|s|}),\qquad(|\arg(s)-\pi|>\theta>0)$$ to (2.5), we
obtain
\begin{equation}
\sum_{\rho}\R\frac{1}{s-\rho}=\frac{1}{2}\log t+O(1).
\end{equation}
We have
$$\aligned \sum_{\rho}\R\frac{1}{s-\rho}=&\sum_{n=1}^{\infty}
\frac{\beta'-\frac{1}{2}}{\left(\beta'-\frac{1}{2}\right)^2+(\gamma'-\gamma_n)^2}+
\frac{\beta'-\frac{1}{2}}{\left(\beta'-\frac{1}{2}\right)^2+(\gamma'+\gamma_n)^2}\\
=&\sum_{n=1}^{\infty}
\frac{\beta'-\frac{1}{2}}{\left(\beta'-\frac{1}{2}\right)^2+(\gamma'-\gamma_n)^2}+O\left(\beta'-\frac{1}{2}\right).
\endaligned$$ Thus, by this and (2.6), we obtain
\begin{equation}
\frac{1}{2}\log\gamma'+O(1)=
O\left(\sum_{n=1}^{\infty}\frac{\beta'-\frac{1}{2}}{(\gamma'-\gamma_n)^2}\right).
\end{equation}
Using (2.4) and Proposition 2.2, we have
$$\sum_{n=1}^{\infty}\frac{\beta'-\frac{1}{2}}{(\gamma'-\gamma_n)^2}=O\left(\left(\beta'-\frac{1}{2}\right)
\frac{(\log\gamma')^2}{\tilde{\delta}^2}\right).$$ We insert this to
(2.7) and then we get
$$\frac{1}{2}\log\gamma'+O(1)=O\left(\left(\beta'-\frac{1}{2}\right)
\frac{(\log\gamma')^2}{\tilde{\delta}^2}\right).$$  This proves
Proposition 2.4.
\end{proof}

Suppose
$$\liminf(\beta'-1/2)\log\gamma'=0.$$
Then this and Proposition 2.4 implies that we have sequences
$\langle\epsilon_k\rangle$, $\langle\rho_{n_k}\rangle$ and
$\langle\rho_k'\rangle$ such that
\begin{equation}
\text{$\rho_{n_k}=1/2+i\gamma_{n_k}$, $\zeta'(\rho_k')=0$,
$|\rho_k'-\rho_{n_k}|<\frac{\epsilon_k}{\log\gamma_{n_k}}$ and
$\epsilon_k\to0\,\,\,(\epsilon_k>0)$.}
\end{equation}
Using (B), we get
$$\frac{1}{\rho_k'-\rho_{n_k}}+O(\log\gamma_{n_k})=0.$$
Thus we obtain that for some $c_1>0$,
$$|\rho_k'-\rho_{n_k}|>\frac{c_1}{\log\gamma_{n_k}}.$$
Note that $c_1$ doesn't depend on $\epsilon_k$'s. By this and (2.8),
we obtain
$$\frac{\epsilon_k}{\log\gamma_{n_k}}>\frac{c_1}{\log\gamma_{n_k}}.$$
But this is a contradiction, since $\epsilon_k\to0$ and $c_1>0$ is a
fixed real number. Thus (A) follows.
\medskip

(A) $\Rightarrow$ (C). Suppose
$$\limsup\frac{|M_n|}{\log\gamma_n}=\infty.$$
We write
$$\frac{\zeta'}{\zeta}(s)=O(\log\gamma_n)+\frac{1}{s-\rho_n}+
\sum_{0<|\gamma_m-\gamma_n|\leqslant1}\frac{1}{s-\rho_m}$$ for
$|s-\rho_n|\leqslant\epsilon/\log\gamma_n$. Let $\epsilon$ be an
arbitrarily small positive real. As in the proof of (C)
$\Rightarrow$ (B), using Proposition 2.2, we obtain
$$\sum_{0<|\gamma_m-\gamma_n|\leqslant1}\frac{1}{s-\rho_m}-\frac{M_n}{i}=O(\log\gamma_n)$$ on
$|s-\rho_n|\leqslant\epsilon/\log\gamma_n$. Then we can see that on
$|s-\rho_n|=\epsilon/\log\gamma_n$ we have
$$\aligned \left|(s-\rho_n)\frac{\zeta'}{\zeta}(s)-(s-\rho_n)
\sum_{0<|\gamma_m-\gamma_n|\leqslant1}\frac{1}{s-\rho_m}\right|
&<\frac{\epsilon}{\log\gamma_n}(M_n+O(\log\gamma_n))\\
&=\left|(s-\rho_n)\sum_{0<|\gamma_m-\gamma_n|\leqslant1}\frac{1}{s-\rho_m}\right|\endaligned$$
for infinitely many $n$'s. We note that
$$(s-\rho_n)\sum_{0<|\gamma_m-\gamma_n|\leqslant1}\frac{1}{(s-\rho_m)}$$ has a zero at $s=\rho_n$.
Thus Rouch\'e's theorem implies that for some $s'$ in
$|s-\rho_n|<\epsilon/\log\gamma_n$,
$$(s'-\rho_n)\frac{\zeta'}{\zeta}(s')=0,$$
i.e., $\zeta'(s')=0$. Therefore $\liminf(\beta'-1/2)\log\gamma'=0$.
Hence (A) implies (C).
\medskip
We have completed the proof of Theorem 4.
\end{proof}

Now we prove Theorem 1 and Corollaries 1, 2.

\begin{proof}[Proof of Theorem 1] Assume that
$$\liminf(\beta'-1/2)\log\gamma'\not=0.$$ Then RH implies that
$\liminf(\beta'-1/2)\log\gamma'$ is positive. Then, by Theorem A, we
obtain $\liminf(\gamma_{n+1}-\gamma_n)\log\gamma_n>0$. Thus, by
Theorem 4, we have (1) $\Rightarrow$ (2).
\medskip

Assume (2) is true. If there exist multiple zeros for $\zeta(s)$, we
immediately get a contradiction from (2). Thus all zeros of
$\zeta(s)$ are simple. Suppose that
$$\liminf(\gamma_{n+1}-\gamma_n)\log\gamma_n=0.$$
Then there exists a sequence of natural numbers $\langle n_k\rangle$
with $n_k<n_{k+1}$ such that
$$(\gamma_{n_k+1}-\gamma_{n_k})\log\gamma_{n_k}\to0$$
as $k\to\infty$. Using (2), we have
$$\aligned\frac{\zeta'}{\zeta}(s)=&\frac{1}{s-\rho_{n_k+1}}+O(\log
t)\\
=&\frac{1}{s-\rho_{n_k}}+O(\log t)\endaligned$$ at
$s=1/2+i(\gamma_{n_k+1}+\gamma_{n_k})/2$. By this, we obtain
$$\frac{1}{s-\rho_{n_k+1}}-\frac{1}{s-\rho_{n_k}}=O(\log\gamma_{n_k})$$
or
$$\frac{1}{\gamma_{n_k+1}-\gamma_{n_k}}=O(\log\gamma_{n_k}).$$
Namely, we have
$$\frac{1}{(\gamma_{n_k+1}-\gamma_{n_k})\log\gamma_{n_k}}=O(1).$$
This is a contradiction, for
$\lim_{k\to\infty}(\gamma_{n_k+1}-\gamma_{n_k})\log\gamma_{n_k}=0$.
Hence, by Theorem 4, we have (2) $\Rightarrow$ (1).
\medskip

Thus Theorem 1 follows.
\end{proof}

\begin{proof}[Proof of Corollary 1] Assume (1). We fix $c>0$. Then, by
Theorem 4 (B), we have
\begin{equation}
\frac{\zeta'}{\zeta}(s)=O\left(\log(|t|+3)\right)
\end{equation}
on $\sigma=1/2+c/\log(|t|+3)$. Using this, it is not hard to see
that Corollary 1 follows from a Phragm\'en-Lindel\"of argument. We
give the detailed proof for convenience. For the following version
of the Phragm\'en-Lindel\"of Theorem, we refer to \cite[p. 138]{Co}.

\begin{thrmpl} Let $G$ be a simply connected
region and let $f$ be an analytic function on $G$. Suppose there is
an analytic function $\varphi:G\to\mathbb{C}$ which never vanishes
and is bounded on $G$. If $M$ is a constant and
$\partial_{\infty}G=A\cup B$ such that

$$\aligned &(a)\text{ for every $a$ in $A$, $\limsup_{s\to a}|f(s)|\leqslant M$;}\\
&(b)\text{ for every $b$ in $B$, and $\eta>0$, $\limsup_{s\to
b}|f(s)||\varphi(s)|^{\eta}\leqslant M$;}\endaligned$$ then
$|f(s)|\leqslant M$ for all $s$ in $G$.
\end{thrmpl}

Here $\partial_{\infty}G=\partial G=\text{ the boundary of }G$ if
$G$ is bounded, $\partial_{\infty}G=\partial G\cup\{\infty\}$ if $G$
is unbounded and the limit superior of $f(s)$ as $s\to a$, is
defined by
$$\limsup_{s\to a}|f(s)|=\lim_{r\to0^+}\sup\{|f(s)|:s\in G\cap
B(a,r)\},$$ where $B(a,r)=\{s\in\mathbb{C}:|s-a|<r\}$. If
$a=\infty$, $B(a,r)$ is the ball in the metric of
$\mathbb{C}_{\infty}=\mathbb{C}\cup\{\infty\}$.

We apply this theorem for proving (1) $\Rightarrow$ (3). We define
sets $G_1$ and $G_2$ by
$$G_1=\{s\in\mathbb{C}:\sigma\geqslant\frac{1}{2}+\frac{c}{\log(|t|+3)}\}\text{ and }G_2=
\{s\in\mathbb{C}:0\leqslant\sigma\leqslant1+c,\,|t|\leqslant1\}.$$
We define $G$ by
$$G=G_1\cap(\mathbb{C}-G_2).$$
We note that $1\not\in G$. We define $f(s)$ by
$$f(s)=\frac{\zeta'(s)}{\zeta(s)\log s}$$
for $s\in G$. Then, $f(s)$ is analytic on the region $G$. We choose
$\varphi(s)=\exp(-\sqrt{s})$. Clearly, the function
$\varphi:G\to\mathbb{C}$ is an analytic function which never
vanishes and is bounded on $G$. By (2.9), there exists $M>0$ such
that we have
\begin{equation}
|f(s)|\leqslant M
\end{equation}
on the boundary $\partial G$.

\begin{claim} We have
$$\frac{\zeta'}{\zeta}(s)=O\left(\log^2(|t|+3)\right)$$
for $s=\sigma+it\in G$.
\end{claim}

\begin{proof}[Proof of Claim] We may suppose $\sigma<2$ and $t\geqslant
10$. Using Proposition 2.1(2), (1.1) and the fact that
$|\sigma-1/2|>c/\log(|t|+3)$ for $s\in G$, we get
$$\aligned\frac{\zeta'}{\zeta}(s)=&\sum_{|t-\gamma|<1}\frac{1}{s-\rho}+O(\log
t)\\
=&O\left(\log t\sum_{|t-\gamma|<1}1\right)+O(\log
t)\\
=&O\left(\log^2t\right).\endaligned$$ Thus, Claim follows.
\end{proof}

Let $\eta>0$. Then, by Claim, we conclude that in $G$, we have
$$\limsup_{s\to\infty}|f(s)||\varphi(s)|^{\eta}=
O\left(\limsup_{r\to\infty}\log
r\exp\left(-\eta\sqrt{r}\cos\frac{\pi}{4}\right)\right)=0.$$ With
this and (2.10), we see that the functions $f$ and $\varphi$ fulfill
the conditions of Phragm\'en-Lindel\"of Theorem. Hence
$|f(s)|\leqslant M$ for $s\in G$. Namely, we obtain
$$\frac{\zeta'}{\zeta}(s)=O(\log|s|)$$
for $s\in G$. In particular, we have
$$\frac{\zeta'}{\zeta}(s)=O(\log t)$$ uniformly for $1/2+c/\log
t\leqslant\sigma\leqslant1/2+c/\log\log t$ $(t\geqslant10)$. From
this, (1.3) and the fact that for $1/2+c/\log
t\leqslant\sigma\leqslant1/2+c/\log\log t$,
$$(\log t)^{2-2\sigma}=O(\log t),$$
we prove Corollary 1.
\end{proof}

\begin{proof}[Proof of Corollary 2] Using (1.4) and Corollary 1, we have
$$\aligned\log\zeta(s)=&\log\zeta(\sigma_1+it)-\int_{\sigma}^{\sigma_1}
\frac{\zeta'(\widetilde{\sigma}+it)}{\zeta(\widetilde{\sigma}+it)}d\widetilde{\sigma}\\
=&O\left(\frac{(\log t)^{2-2\sigma_1}}{\log\log
t}\right)+O\left(\int_{\sigma}^{\sigma_1}(\log t)^{2-2\widetilde{\sigma}}d\widetilde{\sigma}\right)\\
=&O\left(\frac{(\log t)^{2-2\sigma_1}}{\log\log
t}\right)+O\left(\frac{(\log t)^{2-2\sigma}}{\log\log t}\right)\\
=&O\left(\frac{(\log t)^{2-2\sigma}}{\log\log t}\right)\endaligned$$
for $1/2+c/\log t\leqslant\sigma\leqslant\sigma_1$. Thus, we prove
Corollary 2.
\end{proof}

\section{Proof of Theorem 2}

We let $T>2$ and $\sigma=1/2+a/\log T$ for a small $a$. We set
$$a_n=\frac{\gamma_n+\gamma_{n-1}}{2}.$$
Theorem 4 implies that under the assumption of Theorem 2, we obtain
that for $s=\sigma+it$ and $n=2,3,\ldots$,
\begin{equation}
\frac{\zeta'}{\zeta}(s)=\frac{1}{s-\rho_n}+O(\log t)
\end{equation}
for $0\leqslant\sigma-1/2<1/\log\gamma_n$ and $a_n< t\leqslant
a_{n+1}$. We write
$$\aligned
\int_1^T\left|\frac{\zeta'}{\zeta}(\sigma+it)\right|^2dt=&
\sum_{2\leqslant n\leqslant n_1}\int_{a_n}^{a_{n+1}}
\left|\frac{\zeta'}{\zeta}(\sigma+it)\right|^2dt+\int_{a_{n_1+1}}^T\left|\frac{\zeta'}{\zeta}(\sigma+it)\right|^2dt+O(1)\\
=&I+II+O(1),
\endaligned$$
where $a_{n_1+1}< T<a_{n_1+2}$. Using (3.1), it is easy to see that
$$II=\int_{a_{n_1+1}}^T\left|\frac{1}{(\sigma-1/2)+i(t-\gamma_{n_1+1})}
+O(\log T)\right|^2dt=O\left(\frac{\log^2 T}{a}\right)+
o(\log^2T).$$ We insert (3.1) into $I$ and then we get
\begin{equation}\begin{split}
I=&\sum_{2\leqslant n\leqslant n_1}\int_{a_n}^{a_{n+1}}
\left|\frac{1}{(\sigma-1/2)+i(t-\gamma_n)}+O(\log
T)\right|^2dt\\
=&2\sum_{1<\gamma_n<T}\frac{1}{\sigma-1/2}\tan^{-1}\frac{\gamma_{n+1}-
\gamma_n}{2(\sigma-1/2)}+\\
&O\left(\log
T\sum_{1<\gamma_n<T}\log\left(1+\frac{\gamma_{n+1}-\gamma_n}{2(\sigma-1/2)}\right)\right)+
O\left(T\log^2T\right).
\end{split}\end{equation}
Since  $\liminf(\gamma_{n+1}-\gamma_n)\log\gamma_n>0$, there is a
$\beta>0$ such that
\begin{equation}
\frac{\gamma_{n+1}-\gamma_n}{2(\sigma-1/2)}\geqslant\frac{\beta}{2a}
\end{equation}
is large as $a\to0$. We recall Proposition 2.1(1)
\begin{equation}
\sum_{1<\gamma_n<T}1=\frac{T}{2\pi}\log\frac{T}{2\pi e}+O(\log T).
\end{equation}
Using (3.3), (3.4) and the fact that
$\tan^{-1}x=\frac{\pi}{2}+O(1/x)$ as $x\to\infty$, we get
\begin{equation}
\sum_{1<\gamma_n<T}\tan^{-1}\frac{\gamma_{n+1}-
\gamma_n}{2(\sigma-1/2)}=\frac{1}{4}T\log T-\frac{T}{4}\log(2\pi
e)+O(aT\log T)+O(\log T).
\end{equation}
We recall that there exists a $A>0$ such that
\begin{equation}
\#\{n:0<\gamma_n\leqslant T,\gamma_{n+1}-\gamma_n\geqslant
\frac{\lambda}{\log T}\}=O\left(T\log T
e^{-A\lambda^{\frac{1}{2}}(\log\lambda)^{-\frac{1}{4}}}\right),
\end{equation}
uniformly for $\lambda\geqslant2$. For this we refer to \cite{Fu}
and \cite[p. 246]{Ti}. Using (3.4) and (3.6), we can see that
\begin{equation}\begin{split}
\sum_{1<\gamma_n<T}\log\left(1+\frac{\gamma_{n+1}-\gamma_n}{2(\sigma-1/2)}\right)
&=O\left(T\log T\sum_{m=2}^{\infty}\log\left(1+\frac{m}{2a}\right)e^{-Am^{\frac{1}{2}}(\log m)^{-\frac{1}{4}}}\right)\\
&=O\left(T\log T\log\frac{1}{a}\right).
\end{split}\end{equation}
We insert (3.5) and (3.7) into (3.2) and then we obtain
$$I=\frac{1}{2a}T\log^2T-\frac{\log(2\pi e)}{2a}T\log T+O\left(\log\frac{1}{a}(T\log^2T)\right)+O(\log^2 T/a).$$
Using $I$ and $II$, we can get
$$\int_1^T\left|\frac{\zeta'}{\zeta}(\sigma+it)\right|^2dt
=\frac{1}{2a}T\log^2T\left(1-\frac{\log(2\pi e)}{\log T}
+O\left(a\log\frac{1}{a}\right)+O\left(\frac{1}{T}\right) \right).$$
Hence Theorem 2 follows.

\section{Acknowledgment}

I truly thank Professor C. Y. Yildirim for his many valuable
comments and suggestions on this paper.

\end{document}